\magnification=\magstep1
\baselineskip=12pt
\parskip=6pt
\documentstyle{amsppt}
\UseAMSsymbols
\voffset=-3pc
\loadbold
\loadmsbm
\loadeufm
\UseAMSsymbols

\def\cH{\Cal H}
\def\cF{\Cal F}
\def\bC{\Bbb C}

\def\bR{\Bbb R}

\def\var{\varepsilon}
\def\opartial{\overline\partial}
\def\id{\text{id}}
\def\oS{\overline S}
\def\Im{\,\text{Im}\,}
\def\rk{\,\text{rk}\,}
\NoBlackBoxes
\topmatter
\title Geodesics in the space of K\"ahler metrics\endtitle
\author L\'aszl\'o Lempert\footnote{The
first author is grateful to the NSF for partial support
(DMS 0700281) and to the
Universit\'e Pierre et Marie Curie, Paris, where some of this
research was done.\hfill\break} and Liz Vivas\endauthor
\leftheadtext{Space of K\"ahler metrics}
\rightheadtext{L\'aszl\'o Lempert and Liz Vivas}
\address Department of Mathematics, Purdue University, West Lafayette, IN 47907, USA
\endaddress
\subjclassyear{2000}
\subjclass 32Q15, 32W20\endsubjclass
\abstract Let $(X,\omega)$ be a compact K\"ahler manifold. As
discovered in the late 1980s by Mabuchi, the set $\cH_0$ of
K\"ahler forms cohomologous to $\omega$ has the natural structure
of an infinite dimensional Riemannian manifold.
We address the question whether
any two points in $\cH_0$ can be connected by a smooth geodesic, and
show that the answer, in general, is ``no".\endabstract
\endtopmatter
\document

\subhead 1.\ Introduction\endsubhead

Let $X$ be a connected, compact, complex manifold and $\omega_0$ a smooth K\"ahler form on it.
It was discovered by Mabuchi, and rediscovered by Semmes and Donaldson, that
the set $\cH_0$ of smooth K\"ahler forms cohomologous to $\omega_0$, and
the set $\cH$ of smooth, strongly $\omega_0$--plurisubharmonic functions on $X$ carry 
natural infinite dimensional Riemannian manifold structures, see
\cite{M,S,D1}.---A function $u\colon X\to\bR$ is (strongly) $\omega_0$--plurisubharmonic 
if $\omega_0+i\partial\overline\partial u\geq 0$
(resp.~$>0$).---Mabuchi shows that in fact $\cH$ is isometric to the Riemannian product 
$\cH_0\times\bR$, and both he and Donaldson point
out that understanding geodesics in these spaces is important for the study of special 
K\"ahler metrics.
Donaldson then raises the obvious question whether any pair of points in $\cH$ (or $\cH_0$) 
can be connected by a smooth geodesic.
In this paper we give a negative answer:

\proclaim{Theorem 1.1}Suppose $(X,\omega_0)$ is a positive dimensional
connected compact K\"ahler manifold and $h\colon X\to X$ is a holomorphic isometry
with an isolated fixed point such that $h^2=id_X$.
Then there is a K\"ahler form $\omega_1\in\cH_0$ which cannot be connected to $\omega_0$ by a smooth geodesic.
\endproclaim

Concretely, one can take $X$ to be a torus $\bC^m/\Gamma$, $\omega_0$ a translation invariant K\"ahler form, and $h$ induced by reflection
$z\mapsto -z$ in $\bC^m$.

According to Semmes, geodesics in $\cH$ (and therefore in $\cH_0$) are related 
to a 
Monge--Amp\`ere equation as follows, [S].
Let $S=\{s\in\bC\,\colon\, 0<\Im s<1\}$, and $\omega$ the pull back of $\omega_0$ by the projection $\overline S\times X\to X$.
With any smooth curve $[0,1]\ni t\mapsto v_t\in\cH$ associate the smooth function $u(s,x)=v_{\Im s}(x)$, $(s,x)\in\overline S\times X$.
Set $m=\dim X$.
Then $t\mapsto v_t$ is a geodesic if and only if $u$ satisfies
$$
(\omega+i\partial\opartial u)^{m+1}=0.\tag1.1
$$
Since $\omega+i\partial\opartial u$, restricted to fibers $\{s\}\times X$, is positive, (1.1) is equivalent to $\rk \omega+i\partial\opartial
u\equiv m$; and so a smooth geodesic connecting $0, v\in\cH$ gives rise to an $\omega$--plurisubharmonic $u\in C^\infty(\oS\times X)$ solving 
$$
\aligned
 \rk\ \omega+i\partial\opartial u &\equiv m,\\
 u(s+\sigma, x)&=u(s,x)\qquad\text{for }\sigma\in\bR,\ (s,x)\in\oS\times X,\\
 u(s,x)&=\cases 0,&\text{if $\Im s=0$}\\ v(x),&\text{if $\Im s=1$.}\endcases
\endaligned\tag1.2
$$
Therefore Theorem 1.1 follows from the following more precise result:

\proclaim{Theorem 1.2}If $(X,\omega_0)$ and $h$ are as in Theorem 1.1, there is a $v\in\cH$ for which (1.2) admits no real valued solution 
$u\in C^3(\oS\times X)$.
One can choose $v$ so that $h^* v=v$.
\endproclaim

When $m=1$, the $v$ in Theorem 1.2 even form an open subset of the space of $h$--invariant functions in $\cH$, but we do not know if this holds
when $m>1$.

The idea that symmetries help in the analysis of solutions of Monge--Amp\`ere equations is not new.
The first examples of irregularity of certain boundary value problems in $\bC^m$ were constructed by Bedford and Fornaess using symmetries,
see \cite{BF}.
Our approach, based on the study of the so called Monge--Amp\`ere foliation, 
is different from theirs.
The symmetry will be used to identify a leaf of the foliation associated to a $C^3$ solution $u$ of (1.2).
By analyzing the first order behavior of the foliation about this particular leaf we obtain a condition on the Hessian of $u$ at $(1,x_0)$,
where $x_0$ is an isolated fixed point of $h$.
The proof is concluded by finding a boundary value $v$ which is incompatible with this condition.

Studying solutions of the homogeneous Monge--Amp\`ere equation through the 
associated foliation is not new, either.
This approach first appeared in \cite{Be1-2, L1-2}, and still seems to be the only way 
to prove smoothness of the solution.
More recently, in \cite{D2} Donaldson used the foliation method  in 
a variant of the boundary value problem (1.2) to prove,
resp.~disprove, regularity, depending on the boundary data.

Generalized solutions to (1.2) and to rather more general boundary value problems for the homogeneous Monge--Amp\`ere equation (1.1) are known
to exist, see \cite{C}, with complements in \cite{B\l}.

Most of this paper is devoted to the proof of Theorem 1.2. In section 4 we
will discuss the implications of the theorem on the precise regularity
of geodesics in $\cH$.

\medskip\noindent
{\sl Acknowledgement}.
We are grateful to B.~Berndtsson and Z.~B\l ocki for sparking our interest 
in the subject.

\subhead 2.\ The Monge--Amp\`ere foliation\endsubhead

Let $Y$ be an $m+1$ dimensional complex manifold, $\omega$ a real $(1,1)$ form on it, of class $C^1,\ d\omega=0$.
If $\rk\omega\equiv m$, the kernels of $\omega$ form an integrable subbundle of $TY$, and so $Y$ is foliated by Riemann surfaces, whose tangent
spaces are the kernels of $\omega$.
The foliation is of class $C^1$.
If $w$ is a locally defined potential of $\omega$, i.e.~$i\partial\opartial w=\omega$, the section $\partial w$ of $T^{*1,0} Y$ is holomorphic
along any leaf of the foliation (see \cite{BK, Theorem 2.4} for a more general statement).
Applying this with $\omega$ replaced by $\omega+i\partial\opartial u$, we obtain

\proclaim{Proposition 2.1}Suppose $u\in C^3(\oS\times X)$ is a real solution of (1.2).
Then there is a foliation $\cF_u$ on $\oS\times X$, of class $C^1$, whose leaves are Riemann surfaces, tangent to 
$\text{Ker}\,(\omega+i\partial\opartial u)$.
If $w$ is a locally defined potential of $\omega$, then $\partial(w+u)$ is 
holomorphic along the leaves of $\cF_u$.
\endproclaim

$\cF_u$ is called the Monge--Amp\`ere foliation associated with $u$.
We will also need

\proclaim{Proposition 2.2}Suppose $h\colon X\to X$ is a holomorphic map with a fixed point $x_0$.
If $h^2=\id_X$, there are holomorphic local coordinates $z_j$ centered at $x_0$ in which $h$ is expressed as $(z_j)\mapsto (\pm z_j)$.
All signs will be minus if $x_0$ is an isolated fixed point.
\endproclaim

\demo{Proof}This is again not new.
For any compact group of holomorphic transformations local coordinates centered at any fixed point can be found in which the transformations
are linear, see e.g. [BM, p.~19]; and each linear transformation in this compact group 
is diagonalizable.
In our case the eigenvalues of the linearization of $h$ must be $\pm 1$, which means that in an eigenbasis it is indeed given by $(z_j)\mapsto
(\pm z_j)$.
Clearly, the fixed point is isolated only if all the signs are minus.
\enddemo

>From now on we assume $X,\omega_0$, and $h$ are as in Theorem 1.2, and denote by 
$x_0$ an isolated fixed point of $h$.
Let $\tilde h=\id_{\oS}\times h\colon \oS\times X\to\oS\times X$.
We also fix $v\in\cH$ such that $h^* v=v$.

\proclaim{Proposition 2.3}If a real $u\in C^3(\oS\times X)$ solves (1.2), then 
$\tilde h^* u=u$ and $\oS\times \{x_0\}$ is one leaf of the
Monge--Amp\`ere foliation $\cF_u$.
\endproclaim

\demo{Proof}First observe that the solution to (1.2) is unique.
Indeed, if $R=\{s\in\bC\,\colon\, 1 < |s| < e\}$, any solution $u$ of (1.2) defines a solution $U(s,x)=u(i\log s,x)$ of the homogeneous
Monge--Amp\`ere equation on $\overline R\times X$.
The latter being a compact manifold with boundary, uniqueness for the boundary value problem on it is standard. One can argue as follows (see 
\cite{D1, Corollary 7} for a slightly weaker statement). Denote by $\Omega$
the pullback of $\omega_0$ by the projection $\overline R\times X\to X$. The
rank condition implies that the signature of $\Omega+i\partial\bar\partial U$
is constant in $\overline R\times X$; since $\Omega+i\partial\bar\partial U$ is
semipositive at boundary points, it must be semipositive everywhere.
Now suppose $U'$ also solves 
$\text{rk}\, \Omega+i\partial\bar\partial U'\equiv m$, and $U=U'$ on
$\partial R\times X$. Upon adding a constant we can assume $U\ge 0$.
With $\lambda>1$ and $U''=\lambda U'-U$
$$
(\lambda \Omega+i\partial\bar\partial U+i\partial\bar\partial U'')^{m+1}=
\lambda^{m+1}(\Omega+i\partial\bar\partial U')^{m+1}=0.
$$
As $\lambda\Omega+i\partial\bar\partial U$ is positive on the fibers 
$\{s\}\times X$, by [D1, Lemma 6] $U''$ attains its minimum on
$\partial R\times X$. (Donaldson speaks of maximum because he works with
the operator $\bar\partial\partial$ rather than $\partial\bar\partial$.) 
But $U''\ge 0$ on $\partial R\times X$, hence everywhere, and so $U\le\lambda U'$;
letting $\lambda\to 1$ gives $U\le U'$. Reversing the roles of $U$ and $U'$
we see $U= U'$.

Since $u$ and $\tilde h^* u$ solve the same boundary value problem, we conclude that they are indeed equal.
It follows that $\text{Ker}\,(\omega+i\partial\opartial u)$ is invariant under 
$\tilde h_*$, and so is 
Ker$(\omega+i\partial\opartial u)|_{(s,x_0)}\subset
T_{(s,x_0)}(\overline S\times X)$, for any $s\in\overline S$.
Proposition 2.2 shows that the only $\tilde h_*$--invariant complex lines in $T_{(s,x_0)}(\overline S\times X)\approx T_s\overline S\oplus
T_{x_0} X$ are $T_s\overline S\oplus (0)$ and lines contained in $(0)\oplus T_{x_0} X$.
Since $u|(\{0\}\times X)$ is strongly $\omega_0$--plurisubharmonic, Ker$(\omega+i\partial\opartial u)|_{(s,x_0)}$ must agree with $T_s\overline
S\oplus (0)$ when $s=0$.
By continuity, the same must hold for $s$ in a connected neighborhood $S'\subset\overline S$ of $0$, and so $S'\times \{x_0\}$ is contained
in a leaf of $\cF_u$.
By analytic continuation, $\oS\times \{x_0\}$ is itself a leaf, q.e.d.
\enddemo

In the rest of this section, for a real solution 
$u=\tilde h^* u\in C^3(\oS\times X)$ of (1.2) we will study the first 
order behavior of $\cF_u$ about the leaf
$\overline S\times \{x_0\}$.
We fix local coordinates $z_j,\ j=1,\ldots,m$, on a neighborhood $V\subset X$ of $x_0$, 
as in Proposition 2.2; we can choose them so that in
addition $\omega_0|_{x_0}=i\sum dz_j\wedge d\bar z_j|_{x_0}$, and
the local coordinates map $V$ on a convex set in $\bC^m$.
Using the coordinates we identify $V$ with its image in $\bC^m$ and 
$x_0$ with $0\in\bC^m$.
Then $\overline S\times V$ is identified with a subset of $\overline S\times\bC^m$, and occasionally we shall write $z_0$ for the $s$
coordinate of a point $(s,z_1,\ldots ,z_m)\in\oS\times V$.

Exhaust $\oS$ by compact subsets, say by rectangles
$$
S_r=\{s\in \oS\,:\, |\text{Re}\,s|\le r\}, \qquad r\in(0,\infty).
$$
Given $r$, if $a\in V$ is sufficiently close to $x_0=0$, the leaf of
$\cF_u|S_r\times X$ passing through $(0,a)$ is the 
graph of a $C^1$ function $f_a\colon S_r\to V$, holomorphic on 
$\text{int}\,S_r$. For example, $f_0\equiv 0$.
Let the components of $f_a$ be $f_{aj}\colon S_r\to\bC$, $j=1,\ldots,m$, and let 
$f_{a0}(s)=s$.
Write 
$$
\omega|\oS\times V=i\sum^m_{j,k=0}\ \omega_{jk} dz_j\wedge d\bar z_k,\quad
\text{ with }\omega_{jk}=0\text{ if }j\text{ or }k=0.
$$
That the image of $(f_{aj})^m_{j=0}$ is tangent to $\text{Ker}\,(\omega+i\partial\opartial u)$ means
$$
\sum^m_{j=0}\{\omega_{jk}(s,f_a(s))+u_{z_j\bar z_k} (s,f_a(s))\} f'_{aj}(s)
=0,\qquad k=0,1,\ldots ,m.\tag2.1
$$
Further, if $\omega_0|V=i\partial\opartial w_0$ and $w(s,x)=w_0(x)$, then $\partial(w+u)$ is holomorphic along the (interior of the) leaves,
i.e.,
$$
w_{z_j} (s,f_a(s))+u_{z_j}(s,f_a(s)),\qquad j=0,1,\ldots, m,\tag2.2
$$
depend holomorphically on $s\in\,\text{int}\,S_r$.

Replace $a\in\bC^m$ by $ta\in\bC^m$, $t\in [0,1]$, and $f_a$ by
$f_{ta}:S_{r(t)}\to V$, where $\lim_{t\to 0}r(t)=\infty$; then
differentiate (2.1), (2.2) with respect to $t$ at $t=0$.
Note that $\varphi_j=df_{(ta)j}/dt|_{t=0}$ are continuous functions on 
$\oS$, holomorphic on $S$; and $\varphi_0\equiv 0$.
The symmetry of $\omega,u$ implies that odd order derivatives of $\omega_{jk}$ and $u$, with respect to $z_l,\bar z_l,l\geq 1$, all vanish
at $(s,0)$.
Hence differentiating (2.1) gives 
$$
\sum_{l=1}^m\{ u_{z_l\bar z_k s}(s,0)\varphi_l(s)+ 
u_{\bar z_l\bar z_k s}(s,0)\overline{\varphi_l(s)} \}
+\sum^m_{j=1}\ \{\delta_{jk}+u_{z_j\bar z_k}(s,0)\}\ \varphi'_j(s)=0,\tag2.3
$$
$k=1,\ldots,m$; and (2.2) gives that for $j=1,\ldots, m$
$$
 \sum^m_{k=1} \{ w_{0z_j z_k}(0)+u_{z_j z_k} (s,0)\}\varphi_k (s)+\sum^m_{k=1}
\{w_{0z_j\bar z_k} (0)+u_{z_j\bar z_k} (s,0) \}\ \overline{\varphi_k(s)}
$$
are holomorphic.
If $u_{zz},u_{z\bar z}$, and $\varphi$ denote the $m\times m$ matrices 
$(u_{z_j z_k}),(u_{z_j\bar z_k})$, and the column vector
$(\varphi_j)$, $j,k=1,\ldots, m$, then this latter says that
$$
u_{zz} (s,0)\varphi(s)+(I+u_{z\bar z} (s,0))\ \overline{\varphi(s)}=
\psi(s)\tag2.4
$$
is holomorphic.

Now $u(s,0)$ depends only on $\Im s$.
If we let
$$
P=I+u_{z\bar z} (i,0)=I+v_{z\bar z} (0)\quad\text{ and }\quad
Q=u_{zz} (i,0)=v_{zz}(0),\tag2.5
$$
(2.4) implies in light of (1.2)
$$
\psi(s)=\cases \overline{\varphi(s)},&\text{if $\Im s=0$}\\
P\overline{\varphi(s)}+Q\varphi(s),&\text{if $\Im s=1$}.\endcases\tag2.6
$$
On the other hand, restricting (2.3) to real $s$ gives
$$
\varphi'(s)=A\varphi(s)+B\overline{\varphi(s)},\qquad s\in\bR,\tag2.7
$$
where
$$
A=-u_{\bar z zs}(0,0)=-A^*,\qquad B=-u_{\bar z\bar zs}(0,0)=
B^T.\tag2.8
$$
That $A$ is skew adjoint follows from the translational invariance of the self adjoint matrix $u_{\bar z z}(s,0)$, while the symmetry of
$B$ is obvious.
Finally we note
$$
\varphi(0)=a.\tag2.9
$$
The equations (2.6), (2.7), (2.9) describe the first order behavior of $\cF_u$ about the leaf $\oS\times\{x_0\}$.
To summarize:

\proclaim{Proposition 2.4}For any $a\in\bC^m$, the solution of the initial value problem (2.7), (2.9) can be extended to a continuous function
$\varphi\colon\oS\to\bC^m$, holomorphic on $S$, and there is another continuous function $\psi\colon\oS\to\bC^m$, also holomorphic on $S$, that
satisfies (2.6).
\endproclaim

\subhead 3.\ The proof of Theorem 1.2\endsubhead

We consider a real solution $u\in C^3(\oS\times X)$ of the boundary value problem (1.2), with $v\in\cH$ satisfying $h^*v=v$, and choose local
coordinates $z_j$ centered at an isolated fixed point $x_0$ of $h$ as in section 2.
The main point of the proof is

\proclaim{Lemma 3.1}Let $P,Q$ be defined by (2.5).
Assume that
$P$ and $Q$ have a common eigenbasis consisting of real vectors $\xi_1,\ldots,\xi_m\in\bR^m$.
If the eigenvalues of
$$
R=(I+P^2-Q\overline Q)P^{-1}
$$
are simple, then at least one of them is $\geq -2$.
\endproclaim

One can even show that all the eigenvalues are $\ge -2$, but this is irrelevant to
our purposes.
We need an auxiliary result:

\proclaim{Proposition 3.2}Let $A=-A^*$ and $B=B^T$ be $m\times m$ complex matrices.
If the block matrix
$$
M=\pmatrix A&B\\ \overline B&\overline A\endpmatrix\tag3.1
$$
has an eigenvector $\pmatrix x\\ y\endpmatrix$ with eigenvalue $\lambda\in\bC$, then either $\lambda^2\in\bR$ or $|x|=|y|$.
\endproclaim

\demo{Proof}One computes
$$
M^2=\pmatrix C&D\\ -D^*&E\endpmatrix,\qquad\text{with }\quad C^*=C,\ E^*=E.
$$
Denoting by $\langle \ , \ \rangle$ the Euclidean inner product on $\bC^m$, for $\xi,\xi',\eta,\eta'\in\bC^m$ we have
$$
\langle C\xi+D\eta,\xi'\rangle-\langle -D^*\xi+E\eta,\eta'\rangle=\langle\xi,C\xi'+D\eta'\rangle-\langle\eta,-D^*\xi'+E\eta'\rangle.
$$
Now substitute
$$
\pmatrix\xi\\ \eta\endpmatrix=\pmatrix \xi'\\ \eta'\endpmatrix=
\pmatrix x\\ y\endpmatrix,\quad\text{ so that }
\pmatrix C&D\\ -D^*&E\endpmatrix \pmatrix \xi\\ \eta\endpmatrix=
\lambda^2 \pmatrix x\\ y\endpmatrix,
$$
to obtain $\lambda^2 (|x|^2-|y|^2)=\overline\lambda^2 (|x|^2-|y|^2)$, and the claim follows.
\enddemo

\demo{Proof of Lemma 3.1}
The assumptions imply that the $\xi_k$ are eigenvectors of $\overline P$ and $\overline Q$, with eigenvalues conjugate to the eigenvalues of
$P$, resp.~$Q$.
As the eigenvalues of $P$ are real, it follows that $P=\overline P$, $Q$, and 
$\overline Q$ commute among themselves.
Define $A,B$ by (2.8) and let $\pmatrix x\\ y\endpmatrix$ be an eigenvector of the matrix (3.1), with eigenvalue $\lambda\in\bC$.
We apply Proposition 2.4.
With $a=x+\overline y$ it is easy to write down the solution of the initial value problem (2.7), (2.9):\ it is 
$$
\varphi(s)=xe^{\lambda s}+\overline y e^{\overline\lambda s},
$$
first for $s\in\bR$, then by analytic continuation, for all $s\in\oS$.
Upon replacing the eigenvector by a multiple, we can arrange that 
$a=x+\overline y\neq 0$.
>From the first equation in (2.6), $\psi(s)=\overline x e^{\overline\lambda s}+ye^{\lambda s}$.
Substituting in the second, and noting that $\overline s=s-2i$ when Im $s=1$, we obtain
$$
P(\overline x e^{-2i\overline\lambda} e^{\overline\lambda s}+y e^{-2i\lambda} e^{\lambda s})+Q(x e^{\lambda s}+
\overline y e^{\overline\lambda s})=\overline x e^{\overline\lambda s}+ye^{\lambda s}.\tag3.2
$$

If $\lambda\in\bR$, we divide through by $e^{\lambda s}$ to get $P\overline a=e^{2i\lambda}(I-Q)a$.
Combining this equation with its complex conjugate yields
$$
\{(I-Q)(I-\overline Q)-P^2 \} a=0.\tag3.3
$$
Let $P\xi_j=p_j\xi_j$, $Q\xi_j=q_j\xi_j$, with $p_j>0$.
The eigenvalues of $(I-Q)(I-\overline Q)-P^2$ are then $|1-q_j|^2-p_j^2$, 
and by (3.3) one of these numbers, say $|1-q_1|^2-p_1^2$, is
zero.
Hence $|q_1|\leq 1+p_1$, and the corresponding eigenvalue of $R$, 
$(1+p_1^2- |q_1|^2)/p_1\geq -2$, as claimed.

If $\lambda\not\in\bR$, then comparing the coefficients of $e^{\lambda s}$ and $e^{\overline\lambda s}$ in (3.2) yields
$$
Qx=(I-e^{-2i\lambda}P)y,\qquad Q\overline y=(I-e^{-2i\overline\lambda} P)\overline x.
$$
Again, we combine these equations and their complex conjugates to get
$$
Q\overline Q x=(I-e^{-2i\lambda}P) (I-e^{2i\lambda} P)x,\qquad
Q\overline Qy=(I-e^{-2i\lambda}P)(I-e^{2i\lambda}P)y,
$$
or
$$
Rx=(e^{-2i\lambda}+e^{2i\lambda})x,\qquad Ry=(e^{-2i\lambda}+e^{2i\lambda})y.\tag3.4
$$
So one eigenvalue of $R$ is $e^{-2i\lambda}+e^{2i\lambda}$.
If $\lambda$ happens to be imaginary, then $e^{-2i\lambda}+e^{2i\lambda}\ge 2$, as
needed.

It remains to take care of the case when $\lambda$ is neither real nor imaginary.
But in fact this case cannot occur.
Indeed, Proposition 3.2 would imply $|x|=|y|$.
Now $\xi_1,\ldots,\xi_m$ are eigenvectors of $R$, with different eigenvalues; in particular, $e^{-2i\lambda}+e^{2i\lambda}$ is a simple
eigenvalue.
Therefore by (3.4) $x$ and $y$ are proportional to each other, and to a real vector $\xi_j$.
At the price of passing to a multiple, we can assume $x\in\bR^m$, and then $y=\alpha x$, $|\alpha|=1$.
We have
$$
\pmatrix A&B\\ \overline B&\overline A\endpmatrix\ \pmatrix x\\ \alpha x\endpmatrix=\lambda \ \pmatrix x\\ \alpha x\endpmatrix,
$$
whence $Ax+\alpha Bx=\lambda x$ and $Bx+\overline\alpha Ax=\overline\alpha\overline\lambda x$.
Multiplying the latter by $\alpha$ and subtracting from the former we obtain $0=(\lambda-\overline\lambda)x$, or $\lambda\in\bR$, a
contradiction.
With this the proof of the lemma is complete.
\enddemo

It is easy to construct self adjoint, resp.~symmetric, matrices $P>0$ and $Q$ that have a common real eigenbasis, but contrary to the claim of
Lemma 3.1, the eigenvalues of $(I+P^2-Q\overline Q)P^{-1}$ are simple and $<-2$.
For example one can take $P=I$ and $Q$ diagonal, with real eigenvalues $q_1>\ldots >q_m>2$.
Theorem 1.2 would then follow from Lemma 3.1, if we could produce a 
strongly $\omega_0$--plurisubharmonic $v\in C^\infty(X)$ such that $h^* v=v$,
and in the local coordinates constructed,
$v_{z\bar z}(x_0)=P-I$, $v_{zz}(x_0)=Q$.
This can indeed be done:

\proclaim{Lemma 3.3}Let $(X,\omega_0)$ be a K\"ahler manifold with almost
complex structure tensor $J\colon TX\to TX$,
let $h\colon X\to X$ be a holomorphic isometry that fixes $x_0\in X$, and
assume $h^2=\id_X$.
Let $q$ be a real quadratic form on $T_{x_0}X$, invariant under $h_*$.
If its Hermitian part 
$$
 q^{1,1}(\xi)=\bigl(q(\xi)+q(J\xi)\bigr)/2,\qquad \xi\in T_{x_0} X,
$$
satisfies
$q^{1,1}(\xi) + \omega_0(\xi,J\xi)>0$ for nonzero $\xi\in T_{x_0} X$, then 
there is a strongly $\omega_0$--plurisubharmonic $v\in C^\infty(X)$ such 
that $h^* v=v$,
$dv(x_0)=0$, whose Hessian at $x_0$ is $q$.
\endproclaim

The lemma is essentially the same as \cite{D2, Lemma 8}, except for the extra ingredient $h$.
We also had to impose the condition on $q^{1,1}$.
The proof to follow is a slight variation on Donaldson's proof.

\demo{Proof}We choose local coordinates $z_j$ centered at $x_0$ in which $h$ is expressed as $(z_j)\mapsto (\pm z_j)$.
It is clear that in some neighborhood $U$ of $x_0$ there is a $w\in C^\infty(U)$ which has all the properties of $v$, except it is not defined
on all of $X$.
Indeed, if
$$
q=\sum_{j,k} a_{jk} dz_j d\bar z_k+\text{ Re }\sum_{j,k} b_{jk} dz_j dz_k|T_{x_0}X,
$$
then $w=(\sum a_{jk} z_j\bar z_k+\text{ Re }\sum b_{jk} z_j z_k)/2$ will do.
By shrinking $U$ and by scaling the coordinates we can arrange that 
$i\partial\opartial w> (\delta-1)\omega_0$ on $U$, with some $\delta>0$,
and that the coordinates map $U$ on a neighborhood of $\{z\in\bC^n\colon \sum |z_j|^2\leq 1\}$.

Next choose a smooth function $\varphi\colon [-\infty,\infty)\to [0,1]$ such that
$$
\varphi(t)=\cases 0,&\text{if $t>0$}\\ 1,&\text{if $t < -1$}.\endcases
$$
Then $\psi(t)=\varphi(\var\log t)$, $\var>0$, defines a smooth function on $[0,\infty)$, supported on $[0,1]$, such that
$$
\psi(t)=1 \quad\text{ if } t < e^{-1/\var},\quad\text{ and }
\quad t|\psi' (t)|,\ t^2 |\psi'' (t)| < C\var,
$$
with some $C$ independent of $\var$.
Thus the smooth function
$$
v=\cases \psi(|z|^2)w&\text{on $U$}\\ 0&\text{on $X\backslash U$}\endcases
$$
satisfies $dv(x_0)=0$, $h^* v=v$, and its Hessian at $x_0$ is $q$.
Further, on $U$
$$
\omega_0+i\partial\opartial v>\omega_0+i\psi (|z|^2)\partial\opartial w-C'\var\omega_0
>(\delta-C'\var)\omega_0,
$$
with $C'$ independent of $\var$.
If $\var$ is small enough, $v$ therefore has all the required properties.
\enddemo

\subhead 4.\ The smoothness of geodesics\endsubhead

Even without delving into technicalities of infinite dimensional
analysis, it is rather clear that if $X,\omega_0$, and
$v\in\cH$ are as in Theorem 1.2, by this theorem there is no
$C^3$ geodesic connecting $0$ and $v$ in $\cH$. It emerged from a
conversation with Berndtsson that there is not even a $C^2$ geodesic,
and in fact not even a geodesic of Sobolev class $W^{1,2}$ (this space
is the largest in which the geodesic equation can be considered);
but to prove this stronger result, we {\sl will} have to delve 
into technicalities.
Modulo technicalities though, the proof is just standard boot strapping for
ordinary differential equations.

One technicality involves replacing $\cH$ by Banach spaces. From now
on $C^k(X)$ will denote the real Banach space of real functions of
class $C^k$ on $X$. If
$k=2,3,\ldots$ and $l=0,1,\ldots$, let
$$
\cH^k=\{v\in C^k(X): \omega_0+i\partial\opartial v>0\},\qquad
\cF^l=\{\Omega\in C^l_{1,1}(X):\Omega>0\},
$$
open subsets of $C^k(X)$, resp. of the Banach space $ C^l_{1,1}(X)$
of real $(1,1)$--forms of class $C^l$. We also write $C^l_{1,0}(X)$ for
the Banach space of $(1,0)$ forms of class $C^l$. Another technicality is
the notion of vector valued Sobolev spaces. Let $B$ be a Banach space
and $J=[0,1]$. The notions of measurability, integrability, and 
$L^p$-integrability of functions $J\to B$ are rather obvious, and were
codified in [Bo]. The weak derivative of 
$\beta\in L^1_{\text {loc}}(\text{int}\,J;B)$
is a function $\gamma\in L^1_{\text{loc}}(\text{int}\,J;B)$, denoted 
$\dot\beta$, such
that $\int_J\varphi\gamma=-\int_J\dot\varphi\beta$ for all 
$\varphi\in C^\infty(J;\bR)$, $\text{supp}\,\varphi\subset (0,1)$.
The only Sobolev space of $B$ valued functions of interest to us is
$$
W^{1,2}(J;B)=\{\beta\in C(J;B)\, :\,\beta\text{ has a weak derivative in }
L^2(J;B)\}.
$$

Mabuchi's Riemannian metric on $\cH$ extends to a Riemannian metric on
$\cH^2$. The equation for geodesics $\beta:J\to\cH$ involves a function
$$
F:\cF^0\times C^0_{1,0}(X)\to C(X),\qquad 
F(\Omega,f)=\langle f,f\rangle_\Omega\, ,\tag4.1
$$
the pointwise squared length of $f$, measured in the Hermitian metric 
determined by $\Omega$. By [M, (3.1)] or [S, (1.5)] the geodesic equation
is 
$$
\ddot{\beta}(t)=
F\big(\omega_0+i\partial\opartial\beta(t),\partial\dot\beta(t)\big),
\tag4.2
$$
and makes sense for $\beta:J\to\cH^2$.

\proclaim{Proposition 4.1} For $l=0,1,\ldots$, $F$ defines a $C^\infty$ map
$\cF^l\times C^l_{1,0}(X)\to C^l(X)$. Furthermore, given a compact
$K\subset\cF^l$, there is a constant $c$ such that
$$
||F(\Omega,f)||_{C^l(X)}\le c||f||^2_{C^l_{1,0}(X)}\, ,\qquad\text{for }
\Omega\in K,\,f\in C^l_{1,0}(X).
$$
\endproclaim

\demo{Proof} It suffices to prove when the second argument $f$ in
$F(\Omega, f)$ is supported in a coordinate neighborhood $N\subset X$. If
in local coordinates 
$$\gather
\Omega=i\sum\Omega_{jk}\,dz_j\wedge d\bar z_k\quad\text{and}\quad
f=\sum f_j\,dz_j,\quad\text{then}\\
F(\Omega,f)=\sum\Omega^{jk}f_j\bar f_k\quad\text{in } N,
\endgather
$$
and $0$ elsewhere, with $(\Omega^{jk})$ denoting the inverse of the 
matrix $(\Omega_{jk})$.
>From this formula both claims of the Proposition follow.\enddemo

By Proposition 4.1,
$F(\omega_0+i\partial\opartial\beta,\partial\dot\beta)\in 
L^1\big(J;C^{k-2}(X)\big)$ if $\beta \in W^{1,2}(J;\cH^k)$. Hence the 
geodesic equation (4.2) can be
considered for $\beta\in W^{1,2}(J;\cH^k)$ in the weak sense:
we require that the weak derivative of 
$\dot\beta\in L^2\big(J;C^k(X)\big)\subset L^2\big(J;C^{k-2}(X)\big)$ should be
$F(\omega_0+i\partial\opartial\beta,\partial\dot\beta)$.

\proclaim{Theorem 4.2} If $X,\omega_0$, and $v$ are as in Theorem 1.2,
and $k\ge 4$,
then there is no $\beta\in W^{1,2}(J;\cH^k)$ that would satisfy
$\beta(0)=0$, $\beta(1)=v$, and the geodesic equation (4.2) in the
weak sense.\endproclaim

\demo{Proof} Of course, it suffices to prove when $k=4$. We will
show that for any weak solution $\beta\in W^{1,2}(J;\cH^4)$ of
(4.2) the function
$$
w(t,x)=\beta(t)(x),\qquad (t,x)\in J\times X,
$$
is in $C^3(J\times X)$. In general, $w\in C^l(J\times X)$ if (and only if,
but we will not need the reverse implication) 
$\beta\in C^{l-j}\big(J;C^j(X)\big)$ for
$j=0,1,\ldots l$. Indeed, when $l=0$, this follows directly from the
definition of continuity. The case $l>0$ can be reduced to the
case $l=0$ by observing that $w\in C^l(J\times X)$ if
$\xi_1\ldots\xi_l w\in C(J\times X)$ for smooth vector fields
$\xi_i$, each tangential either to the fibers $\{t\}\times X$ or to the
fibers $J\times\{x\}$.

By the definition of $W^{1,2}(J;\cH^4)$,
$$
\beta\in C\big(J;C^3(X)\big).\tag4.3
$$
Further, by Proposition 4.1 the map 
$t\mapsto 
F\big(\omega_0+i\partial\opartial\beta(t),\partial\dot\beta(t)\big)\in C^2(X)$
is integrable. Hence by (4.2) $t\mapsto \dot\beta(t)\in C^2(X)$ is
continuous, or rather has a representative, continuous on $J$ (namely the
one obtained by integrating $\ddot{\beta}$). In other words,
$$
\beta\in C^1\big(J;C^2(X)\big).\tag4.4
$$
Feeding this and (4.3) back into (4.2), by Proposition 4.1 we obtain
$\ddot{\beta}= F(\omega_0+i\partial\opartial\beta,\partial\dot\beta)
\in C\big(J;C^1(X)\big)$, or
$$
\beta\in C^2\big(J;C^1(X)\big).\tag4.5
$$
Finally, feeding this and (4.4) back into into (4.2) once more gives
$\ddot\beta\in C^1\big(J;C(X)\big)$, whence $\beta\in C^3\big(J;C(X)\big)$.
This latter, together with (4.3), (4.4), and (4.5), implies
$w\in C^3(J\times X)$, and therefore $u(s,x)=w(\text{Im}\, s,x)$, 
$(s,x)\in\oS\times X$, is in $C^3(\oS\times X)$. However, (4.3-5) also
imply that the two sides of (4.2), viewed as elements of $C(X)$, agree
for every $t\in J$, not just weakly. By the computations in [S] or [D1], this
means that $\text{rk}\,\omega_0+i\partial\opartial u\equiv m$, and of course
the other equations in (1.2) are also satisfied. But if $v\in\cH$ is suitably
chosen, by Theorem 1.2 the solution of (1.2) cannot be in $C^3(\oS\times X)$,
and consequently $\beta$ cannot be in $W^{1,2}(J,\cH^4)$.\enddemo

\NoBlackBoxes

\enddocument
\bye